\documentclass[12pt,reqno]{amsart}

\usepackage{latexsym}
\usepackage[dvips]{graphics,epsfig}
\newtheorem{theorem}{Theorem}
\newtheorem{lemma}{Lemma}
\newtheorem{proposition}{Proposition}
\newtheorem{remark}{Remark}

\newtheorem{corollary}{Corollary}

\newcommand{\beq}{\begin{equation}}
\newcommand{\eeq}{\end{equation}}
\newcommand{\beqna}{\begin{eqnarray*}}
\newcommand{\eeqna}{\end{eqnarray*}}
\newcommand{\beqn}{\begin{equation*}}
\newcommand{\eeqn}{\end{equation*}}
\newcommand{\bp}{\begin{proof}}
\newcommand{\ep}{\end{proof}}
\newcommand{\bprop}{\begin{proposition}}
\newcommand{\eprop}{\end{proposition}}
\newcommand{\bt}{\begin{theorem}}
\newcommand{\et}{\end{theorem}}
\newcommand{\bex}{\begin{Example}}
\newcommand{\eex}{\end{Example}}
\newcommand{\bc}{\begin{corollary}}
\newcommand{\ec}{\end{corollary}}
\newcommand{\bl}{\begin{lemma}}
\newcommand{\el}{\end{lemma}}
\newcommand{\br}{\begin{remark}}
\newcommand{\er}{\end{remark}}

\newcommand{\R}{{\mathbb R}}

\begin{document}

\title
[An extension of polynomial integrability to dual quermassintegrals]
{An extension of polynomial integrability to dual quermassintegrals}

\author{Vlad Yaskin}\thanks{The   author is supported in part by NSERC}

\address{Department of Mathematical and Statistical Sciences, University of Alberta, Ed-
monton, Alberta T6G 2G1, Canada} \email{yaskin@ualberta.ca}

\subjclass[2010]{Primary 52A20, 33C55}

\keywords{Convex bodies, dual quermassintegrals, spherical harmonics}

\begin{abstract}
A body $K$ is called polynomially integrable if its parallel section function $V_{n-1}(K\cap\{\xi^\perp+t\xi\})$ is a polynomial of $t$ (on its support) for every $\xi$. A complete characterization of such bodies was given recently.
Here we obtain a generalization of these results in the setting of dual quermassintegrals. We also address the associated smoothness issues.
\end{abstract}

\maketitle

\section{Itroduction}

Let $K$ be a convex body in $\mathbb R^n$. The {\it parallel section function} of $K$ in  the direction $\xi\in S^{n-1}$ is defined by $$A_{K,\xi} (t)= V_{n-1}(K\cap\{\xi^\perp+t\xi\}),$$
where $\xi^\perp = \{x\in \mathbb R^n: \langle x,\xi\rangle = 0\}$ and $V_{n-1}$ is the $(n-1)$-dimensional volume (Lebesgue measure).

 Following Agranovsky   \cite{A}, we say that a body  $K$ is   {\it  polynomially integrable} if $$A_{K,\xi}(t)   =\sum_{k=0}^N a_k(\xi)\ t^k,$$
  for some integer $N$, all $\xi\in S^{n-1}$ and all $t$ for which
  the set $K\cap\{\xi^\perp+t\xi\}$ is non-empty. Here, $a_k$ are functions on the sphere.

It is not difficult to see that ellipsoids are polynomially integrable in odd dimensions. Agranovsky asked whether there are other polynomially integrable bodies. This question was answered in \cite{A} and \cite{KMY}. In even dimensions there are no polynomially integrable bodies, while in odd dimensions the only polynomially integrable bodies are ellipsoids. Note that these results were obtained under the $C^\infty$ smoothness condition. Below we will address this issue.

Agranovsky's question is motivated by a problem about algebraically integrable domains that goes back to Newton \cite{N}.  The problem is to describe bounded domains with smooth boundaries for which the volume cut off by a plane
from such a domain   depends algebraically on the plane. The   case of convex bodies in $\mathbb R^2$  was settled by Newton, who showed that such bodies do not exist. Arnold (see problems 1987-14, 1988-13, and
1990-27 in \cite{Arnold}) asked for extensions of Newton's result to other
dimensions and non-convex domains.
In even dimensions Arnold's question  was answered   by Vassiliev  \cite{V}.

In this paper we   deal with a generalization of the result in \cite{KMY} to dual quermassintegrals. Before we introduce the latter, let us recall some facts from the classical Brunn-Minkowski theory. In the heart of this theory  is  the behavior of volume under the Minkowski addition. Intrinsic volumes are important concepts in this theory and arise as coefficients in the Steiner formula. If $K$ is a convex body in $\mathbb R^n$ and $B_2^n$ is the unit Euclidean ball, then $$V_n(K+\epsilon B_2^n) = \sum_{i=0}^n{ n\choose i} W_i(K) \epsilon^i=\sum_{i=0}^n\kappa_{n-i} V_i(K) \epsilon^{n-i},$$
where the addition is the Minkowski addition and  $\kappa_i$ is the volume of the unit Euclidean ball in $\mathbb R^i$. The coefficients $ W_i(K)$ and $V_i(K)$ are known as quermassintegrals and intrinsic volumes, correspondingly. Clearly, $V_i(K)$ is the same as  $W_{n-i}(K)$, up to a constant factor. The standard reference for the Brunn-Minkowski theory is the book of Schneider \cite{Sch}.

The dual Brunn-Minkowski theory was initiated by Lutwak \cite{L1}. In this theory the Minkowski addition of convex bodies is replaced by the radial addition of star bodies.
It served as a foundation for many developments in modern convex geometry, especially in problems related to sections of convex or star bodies.
One of the famous examples is the Busemann-Petty problem, where one of the main ingredients were Lutwak's intersection bodies introduced within the dual Brunn-Minkowski theory; see \cite{L2}. For the history of this problem the reader is referred to \cite{K}.

It turns out that there are many analogies between the Brunn-Minkowski theory and its dual counterpart. For example, there is a version of the Steiner formula in the dual setting.
Let $K$ be a star body in $\mathbb R^n$  and let $B_2^n$ be the unit Euclidean ball.
Consider the star body $K\tilde{+} \epsilon B_2^n$, where $\tilde{+}$ is the radial addition.
Then
$$V_n(K\tilde{+} \epsilon B_2^n) = \sum_{i=0}^n{ n\choose i} \tilde W_i(K) \epsilon^{ i}= \sum_{i=0}^n{ n\choose i} \tilde V_i(K) \epsilon^{n-i},$$
where the coefficients $ \tilde W_i(K) $ and  $\tilde V_i(K)$ are called the dual quermassintegrals and dual   volumes, correspondingly. Note that they depend on the choice of the origin.
Explicitly, they can be written as follows:
$$\tilde V_i (K) =  \tilde W_{n-i}(K)= \frac{1}{n} \int_{S^{n-1}} \rho_K^i(\theta) \, d\theta,$$
where the integration is with respect to the spherical Lebesgue measure and $\rho_K$ is the radial function of $K$.

Analogous to the Kubota integral recursion formula, which gives the $i$-th intrinsic volume as the average volume of projections of the body onto $i$-dimensional subspaces,  is its dual version, which some authors call the dual Kubota integral recursion formula,
 \begin{equation}\label{dual_Kubota} \tilde V_i(K) = \frac{\kappa_{n}}{\kappa_i}\int_{G(n,i)} V_i (K\cap   H )  \, dH,
 \end{equation}
 where the integration is  with respect to the Haar probability measure on the Grassmanian $G(n,i)$. This formula, as well as a brief introduction into the dual Brunn-Minkowski theory, can be found in \cite{Ga3}.

 Dual quermassintegrals are particular cases of dual mixed volumes. For some results about applications and characterizations of dual mixed volumes, the reader is referred to \cite{DGP}, \cite{JV}, \cite{Kl}, \cite{M}.

Now we are ready to state our results. Let $K$ be a convex body in $\mathbb R^n$ that contains the origin in its interior and let $m$ be an integer,  $1\le m\le n-1$. Consider the function $A_{K,m, \xi} (t)$ that gives the $m$th dual   volume of the section $K\cap \{\xi^\perp + t\xi\}$  with respect to the point $t\xi$, i.e.,
\begin{align}\label{int_formula}
	A_{K,m, \xi} (t) = \tilde V_m( K\cap \{\xi^\perp + t\xi\} ) & = \frac{1}{n-1} \int_{S^{n-1}\cap  \xi^\perp} \rho_{K-t\xi}^m(\theta)\, d\theta  \notag \\
	&  = \frac{\kappa_{n-1}}{\kappa_m}\int_{G(\xi^\perp, m)} V_m ((K-t\xi) \cap   H )  \, dH,
	\end{align}
where $G( \xi^\perp, m)$ denotes the Grassmanian of $m$-dimensional subspaces of $\xi^\perp$. 
Note that $t\xi$ is not an interior point of the section $ K\cap \{\xi^\perp + t\xi\}$ when $t\ge  \rho_K(\xi)$. The latter is the radial function in the direction $\xi$; see  the next section for its definition. One can also use  representation (\ref{int_formula}) as a definition of $A_{K,m, \xi}(t)$ for all $t$, but we do not really need this.

We say that $K$   has a polynomial $m$th dual section function if there is a positive integer $N$  such that
$$A_{K,m, \xi} (t) = \sum_{k=0}^N a_k(\xi) t^k$$
for all $\rho_K(-\xi)\le t\le \rho_K(\xi) $, where $a_k$ are continuous functions on the sphere.

In this paper we show  that for even $m$ the only convex bodies with this property are  ellipsoids, while  for   odd $m$  such bodies do not exist.
Note that the results of \cite{KMY} can be obtained as a particular case of our results when $m=n-1$. Since we do  not impose any smoothness assumptions on the bodies, it means that the $C^\infty$ assumption used in \cite{KMY} can be dropped. The proof here is based on spherical harmonics, instead of the Fourier transform as in \cite{KMY}.

\section{Definitions and auxiliary facts}

Let $K$ be a compact set in $\mathbb R^n$ that is star-shaped about the origin $o$; the latter means that for every point $x\in K$ the closed segment $[o,x]$ is also contained in $K$. The {\it radial function} of $K$ is defined by \begin{align*}
\rho_K(\xi) = \max\{a \ge 0 :a \xi \in K\}, \quad \xi\in S^{n-1}.
\end{align*}
 We say that $K$ is a {\it star body} if it is compact, star-shaped about the origin and its radial function  is a positive continuous function on $S^{n-1}$.

 The {\it Minkowski functional} of $K$ is defined by
\begin{align*}
\|x\|_K=\min \{a\ge 0: x \in aK \}, \qquad x\in \mathbb R^n.
\end{align*}
Clearly, $\|\cdot\|_K$ is a 1-homogeneous function on $\mathbb R^n$ and  $ \|\xi\|_K=\rho_K(\xi) ^{-1}$ for $\xi\in S^{n-1}$.

If $K$ and $L$ are star bodies, and $\alpha$ and $\beta$ are positive numbers, then the radial sum $\alpha K\tilde{+} \beta L$, which we saw in the introduction, is defined as the star body with the radial function
$$\rho_{\alpha K\tilde{+} \beta L} = \alpha \rho_{K}+\beta \rho_{L}.$$

Let  $h$ be an integrable
function on $\R$ that is $C^\infty$-smooth in a neighborhood of the
origin. If $q\in \mathbb C$, $-1<\Re q<0$, then the {\it fractional derivative}  of the function
$h$ of order $q$ at zero is defined by
$$h^{(q)}(0)=\frac{1}{\Gamma(-q)} \int_0^\infty t^{-1-q}  h(t)  dt.$$
For other values of $q$ it is defined using   analytic continuation (with respect to $q$).
In particular, if $q$ is not an integer and $-1< \Re q <m$ for some
integer $m$, then

$$h^{(q)}(0)=\frac{1}{\Gamma(-q)} \int_0^1 t^{-1-q} \Big(h(t)-h(0)-\cdots -h^{(m-1)}(0) \frac{t^{m-1}}{(m-1)!}\Big) dt$$    $$ +\frac{1}{\Gamma(-q)}
\int_1^\infty t^{-1-q} h(t)
dt+\frac{1}{\Gamma(-q)}\sum_{k=0}^{m-1}\frac{h^{(k)}(0)}{k!(k-q)};$$
see e.g., \cite{K}.

If  $k\ge 0$ is an integer, then passing to the limit in the previous expression we get  ordinary derivatives (up to a sign),
$$h^{(k)}(0) = (-1)^k \frac{d^k}{dt^k} h(t)\Big|_{t=0}.$$

In this paper we will extensively use spherical harmonics; see \cite{Gr} for a comprehensive introduction into this subject.  Recall that a {\it spherical harmonic}  of degree $m$ and dimension $n$ is   the restriction to $S^{n-1}$ of a harmonic  homogeneous polynomial of degree $m$ on $\mathbb R^n$. Given $f\in L^2(S^{n-1})$, one can associate with it its spherical harmonic expansion
$\sum_{m=0}^\infty H_m$, where $H_m$ is a spherical harmonic of order $m$.

The following is the well-known Funk-Hecke theorem. If $\Phi$ is an integrable function on $[-1,1]$ and $H_m$ is a spherical harmonic of order $m$, then
$$\int_{S^{n-1}} \Phi(\langle \xi, \theta\rangle) H_m(\theta)\, d\theta = \lambda_{m}(\Phi) H_m(\xi),$$
where
$$ \lambda_{m}(\Phi) = \int_{-1}^1 \Phi(t) P_m^n(t)  (1-t^2)^{(n-3)/2} \, dt$$
and $P_m^n$ is the Legendre polynomial of degree $m$ and dimension $n$.

\section{Main Results}

\begin{lemma}\label{Lem}
	Let $E$ be an ellipsoid that contains the origin in its interior. Then $E$ has a polynomial $m$th dual section function for every even integer $m$, $2\le m \le n-1$.
	
\end{lemma}

\bp
By formula (\ref{int_formula}), we have
\begin{align*}	A_{E,m, \xi} (t)& =   \frac{\kappa_{n-1}}{\kappa_m}\int_{G(\xi^\perp, m)} V_m ((E-t\xi) \cap   H )  \, dH \\
&=\frac{\kappa_{n-1}}{\kappa_m}\int_{G(\xi^\perp, m)} V_m (E\cap (H\vee\xi) \cap \{\xi^\perp +t\xi\} )  \, dH,
\end{align*}
 where $H\vee \xi= \mbox{span} (H,\xi)$.

We  now use the fact that $E\cap (H\vee\xi)$ is an ellipsoid, and therefore polynomially integrable in $ H\vee \xi$ if $m$ is even (since $ H\vee\xi$ is $(m+1)$-dimensional). Thus

$$ V_m (E\cap  (H\vee\xi) \cap \{\xi^\perp +t\xi\}  ) = \sum_{k=0}^{m}  a_k(H, \xi) t^k,$$ for all $\rho_K(-\xi)\le t \le  \rho_K(\xi)$.

Note that $a_k(H, \xi)$ are continuous functions of $H$ and $\xi$, and thus  denoting $$a_k(\xi) = \frac{\kappa_{n-1}}{\kappa_m}\int_{G(\xi^\perp,m)}a_k(H, \xi) \, dH,$$ we get
$$A_{E,m, \xi} (t) = \sum_{k=0}^{m}  a_k(\xi) t^k,$$ for all $\rho_K(-\xi)\le t \le  \rho_K(\xi)$.

\ep
	
\begin{theorem}
There are no   convex bodies  with polynomial $m$th dual section functions if $m$ is odd.
\end{theorem}

\bp
Let $K$ be a convex body in $\mathbb R^n$ containing the origin as an interior point.
For a fixed $t$, $|t|< \min_{\xi\in S^{n-1}}\rho_K(\xi)$,  consider the integral
$$\int_{S^{n-1}}  A_{K,m, \xi} (t) \, d\xi  = \frac{\kappa_{n-1}}{\kappa_m} \int_{S^{n-1}}   \int_{G(\xi^\perp, m)} V_m (K\cap  (H\vee\xi) \cap \{\xi^\perp +t\xi\} )  \, dH \, d\xi,$$
where $ H\vee \xi= \mbox{span} (H,\xi)$, as before.

Observe that
\begin{align*} &V_m (K\cap (H\vee\xi) \cap \{\xi^\perp + t\xi\} ) = \lim_{\epsilon\to 0^+} \frac{\epsilon}{2}  \int_{K\cap (H\vee\xi)} |\langle x,\xi\rangle - t |^{-1+\epsilon}\, dx \\
& = \lim_{\epsilon\to 0^+} \frac{\epsilon}{2}  \int_{S^{n-1}\cap (H\vee\xi) } \int_0^{\rho_K(\theta)}r^m |r \langle \theta,\xi\rangle - t |^{-1+\epsilon}\,  dr\,  d\theta \\
&= \lim_{\epsilon\to 0^+} \frac{\epsilon}{2}  \int_{S^{n-1}\cap H } \int_{-1}^1 (1-z^2)^{(m-2)/2} \int_0^{\rho_K(z\xi + \sqrt{1-z^2} \eta)}r^m |r z - t |^{-1+\epsilon} \, dr\,  dz\,   d\eta .
\end{align*}

Note that for a function $f$, integrable on the corresponding sub-spheres, one has
$$ \int_{G( \xi^\perp, m)}  \int_{S^{n-1}\cap H } f(\eta)\, d\eta\, dH =
\frac{m \kappa_m}{(n-1) \kappa_{n-1}}  \int_{S^{n-1}\cap   \xi^\perp} f(\eta)\, d\eta.$$
Therefore,
\begin{align*} & \int_{S^{n-1}}  A_{K,m, \xi} (t) \, d\xi \\
&  = \lim_{\epsilon\to 0^+} \frac{\epsilon\, m}{2(n-1)}  \int_{S^{n-1}}  \int_{S^{n-1}\cap   \xi^\perp} \int_{-1}^1 (1-z^2)^{(m-2)/2} \int_0^{\rho_K(z\xi + \sqrt{1-z^2} \eta)}r^m |r z - t |^{-1+\epsilon}\, dr\,  dz\,   d\eta\,  d\xi \\
& = \lim_{\epsilon\to 0^+}\frac{\epsilon\, m}{2(n-1)} \int_{S^{n-1}}  \int_{S^{n-1}   }   (1-\langle \theta,\xi\rangle ^2)^{(m-n+1)/2} \int_0^{\rho_K(\theta )} r^m |r \langle \theta,\xi\rangle - t |^{-1+\epsilon}\,  dr\,  d\theta\,   d\xi\\
& = \lim_{\epsilon\to 0^+} \frac{\epsilon\, m}{2(n-1)} \int_{S^{n-1}}  \int_0^{\rho_K(\theta )} r^m \int_{S^{n-1}   }   (1-\langle \theta,\xi\rangle ^2)^{(m-n+1)/2} |r \langle \theta,\xi\rangle - t |^{-1+\epsilon}\,  d\xi\,   dr \, d\theta   .
\end{align*}
Let us now compute
\begin{align*}    \lim_{\epsilon\to 0^+}& \frac{\epsilon}{2}  \int_{S^{n-1}   }   (1-\langle \theta,\xi\rangle ^2)^{(m-n+1)/2} |r \langle \theta,\xi\rangle - t |^{-1+\epsilon}\,  d\xi    \\
& =  \lim_{\epsilon\to 0^+} \frac{\epsilon}{2} r^{-1+\epsilon}  \int_{S^{n-1} \cap \theta^\perp   }  \int_{-1}^1  (1-z^2)^{(m-2)/2} \left|  z - \frac{t}{r} \right|^{-1+\epsilon}\,  dz\,  d\eta    \\
 & =  \left\{\begin{array}{ll}  r^{-1 }  |S^{n-2}|    \left(1-\frac{t^2}{r^2}\right)^{(m-2)/2} , & \mbox{ if } t\le r  ,\\
 0, & \mbox{ if } t>r .\end{array}
 	\right.
 \end{align*}

Thus,
\begin{align*}
\int_{S^{n-1}}  A_{K,m, \xi} (t) \, d\xi& =  m\kappa_{n-1} \int_{S^{n-1}}  \int_t^{\rho_K(\theta )}  r^{m-1 }     \left(1-\frac{t^2}{r^2}\right)^{(m-2)/2} \,  dr \,  d\theta \\
&=  m\kappa_{n-1}\int_{S^{n-1}}  \int_t^{\rho_K(\theta )}  r    \left(r^2 - t^2 \right)^{(m-2)/2}  \,   dr \,  d\theta \\
& =  \kappa_{n-1}  \int_{S^{n-1}}      \left(\rho_K^2 (\theta ) - t^2 \right)^{m/2}  \,    d\theta .
\end{align*}

If $K$ had a polynomial $m$th dual section function, then $\int_{S^{n-1}}  A_{K,m, \xi} (t) \, d\xi$ would be a polynomial of $t$, for $|t|< \min_{\xi\in S^{n-1}}\rho_K(\xi)$. But, clearly, $ \int_{S^{n-1}}      \left(\rho_K^2 (\theta ) - t^2 \right)^{m/2}  \,    d\theta$ is not a polynomial, since $m$ is odd.

\ep

\begin{theorem}
For   even $m$, the only   convex bodies with a polynomial $m$th dual section function  are ellipsoids.
\end{theorem}

\bp
Let $K$ be convex body in $\mathbb R^n$ that contains the origin in its interior. Assume that $K$ has a polynomial $m$th dual section function for some even $m$.

For a fixed $\xi$, let us compute fractional derivatives of $ A_{K,m, \xi} (t) $ at zero.
Assume that $-1<\Re q<0$. Then
$$ A_{K,m, \xi}^{(q)} (0) = \frac{1}{\Gamma(-q)} \int_0^\infty t^{-1-q}  A_{K,m, \xi} (t)\, dt$$
$$= \frac{1}{2\Gamma(-q)} \int_{-\infty}^\infty |t|^{-1-q} (1+\mbox{sgn}\, t)  A_{K,\xi,m} (t) \, dt$$
\begin{align*} &=  \frac{\kappa_{n-1}}{\kappa_m}  \frac{1}{2\Gamma(-q)}  \int_{G(\xi^\perp, m)}   \int_{-\infty}^\infty |t|^{-1-q} (1+\mbox{sgn}\, t)   V_m (K\cap (H\vee\xi) \cap \{\xi^\perp +t\xi\} ) \, dt  \, dH \\
&=  \frac{\kappa_{n-1}}{\kappa_m}  \frac{1}{2\Gamma(-q)}  \int_{G(\xi^\perp, m)} \int_{ K\cap  (H\vee\xi)} |\langle x,\xi\rangle|^{-1-q} (1 + \mbox{sgn}\, \langle x,\xi\rangle) \,  dx  \, dH \\
&=  \frac{\kappa_{n-1}}{\kappa_m}  \frac{1}{2\Gamma(-q)}  \int_{G(\xi^\perp, m)}  \int_{ S^{n-1}\cap (H\vee\xi)} \int_0^{\rho_K(\theta)} r^{m-1-q} dr |\langle \theta,\xi\rangle|^{-1-q}   (1 + \mbox{sgn}\, \langle \theta,\xi\rangle)\,  d\theta  \, dH \\
&=  \frac{\kappa_{n-1}}{\kappa_m}  \frac{1}{2\Gamma(-q) (m-q)}  \int_{G(\xi^\perp, m)}  \int_{ S^{n-1}\cap  (H\vee\xi)}  \rho_K^{m-q}(\theta)    |\langle \theta,\xi\rangle|^{-1-q}  (1 + \mbox{sgn}\, \langle \theta,\xi\rangle) \,  d\theta  \, dH \\
&=  \frac{\kappa_{n-1}}{\kappa_m}  \frac{1}{2\Gamma(-q) (m-q)}  \int_{G(\xi^\perp, m)}  \int_{ S^{n-1}\cap  H}  \int_{-1}^1 \rho_K^{m-q}(z\xi + \sqrt{1-z^2} \eta)    |z|^{-1-q} (1 + \mbox{sgn}\, z) \\
&\qquad \hspace{9cm} \times (1-z^2)^{(m-2)/2} \, dz \, d\eta \, dH \\
&=     \frac{m}{2\Gamma(-q) (m-q)(n-1)}    \int_{ S^{n-1}\cap    \xi^\perp}  \int_{-1}^1 \rho_K^{m-q}(z\xi + \sqrt{1-z^2} \eta)    |z|^{-1-q} (1 + \mbox{sgn}\, z)  \\
 &\qquad \hspace{9cm} \times (1-z^2)^{(m-2)/2} \, dz\,  d\eta   \\
&=     \frac{m}{2\Gamma(-q) (m-q)(n-1)}    \int_{ S^{n-1} }    \rho_K^{m-q}(\theta)    |\langle \theta,\xi\rangle|^{-1-q}  (1 + \mbox{sgn}\, \langle \theta,\xi\rangle)   \\
 &\qquad \hspace{9cm} \times  (1-\langle \theta,\xi\rangle^2)^{(m-n+1)/2} \,  d\theta  .
\end{align*}

For $-1<\Re q <0$ and  a function $f\in C(S^{n-1})$ we will be interested in the operator $I_q$ defined by
$$I_qf (\xi) =   \frac{1}{2\Gamma(-q) }    \int_{ S^{n-1} }       |\langle \theta,\xi\rangle|^{-1-q}  (1 + \mbox{sgn}\, \langle \theta,\xi\rangle)    (1-\langle \theta,\xi\rangle^2)^{(m-n+1)/2}  f(\theta) d\theta  . $$
  If $H_l$ is a spherical harmonic of degree $l$, then by Funk-Hecke's theorem, $$I_q H_l = \lambda_l(q) H_l,$$ where
$$\lambda_l(q) =   \frac{1}{ \Gamma(-q) }   \int_{0}^1 P^n_l(z)     z ^{-1-q}    (1-z^2)^{(m-2)/2} dz  , $$
where $P^n_l$ is the Legendre polynomial of degree $l$ and dimension $n$.

Such multipliers were studied before, e.g., by Goodey \cite{Goodey}. But he was interested in the case when $-1-q$ is a non-negative integer, while our main interest is when $-1-q$ is a negative integer.

Note that $ \lambda_l(q)$ is the fractional derivative of order $q$  (with $-1<\Re q <0$) at zero of  the function  $P^n_l(z)  (1-z^2)^{(m-2)/2}$. Thus $ \lambda_l(q)$ can be extended to all  $q$ using analytic continuation.
So, if $q$ is equal to an integer $k\ge 0$, then
$$\lambda_l(k) = (-1)^k \frac{d^k}{dz^k} \left( P^n_l(z)      (1-z^2)^{(m-2)/2} \right)\Big|_{z=0}.$$

Note that $P^n_l(z)      (1-z^2)^{(m-2)/2}$ is a polynomial of degree $l+m-2$ and thus $\lambda_l(k) = 0 $ if $m<k-l+2$. In addition, $\lambda_l(k) = 0 $ if $l$ and $k$ are of different parity. Let us now show that $\lambda_l(k) \ne 0 $ for all $m\ge k-l+2$ if $k-l$ is even.

We have
\begin{align*} \lambda_l(k) & = (-1)^k \frac{d^k}{dz^k} \left( P^n_l(z)      (1-z^2)^{(m-2)/2} \right)\Big|_{z=0}\\
& = (-1)^k \sum_{j=0}^k \frac{d^{k-j}}{dz^{k-j}}  P^n_l(z)     \Big|_{z=0}  \,  \frac{d^{j}}{dz^{j}} (1-z^2)^{(m-2)/2}    \Big|_{z=0}  .
\end{align*}
Since
$$  (1-z^2)^{(m-2)/2}    = \sum_{i=0}^{(m-2)/2}   (-1)^i {{ (m-2)/2} \choose i} z^{2i},$$
it follows that
$$\frac{d^{j}}{dz^{j}} (1-z^2)^{(m-2)/2}     \Big|_{z=0} =
 (-1)^{j/2} {{ (m-2)/2} \choose  j/2} j! .$$
 if $j$ is even and $j\le m-2$, and this derivative is zero if $j$ is odd.

Also, by \cite[Lemma 3.3.9]{Gr}, we have
$$ \frac{d^{k-j}}{dz^{k-j}}   P_l^n (z)  = c_{n,l,j} P^{n+2k-2j}_{l-k+j}(z),$$
 if $k-j\le l$ and it is zero otherwise. Here $c_{n,l,j}$ is a positive constant.

Now use   \cite[Lemma 3.3.8]{Gr},

$$ P^{n+2k-2j}_{l-k+j}(0) = (-1)^{(l-k+j)/2} b_{n,l,j},$$
if $l-k+j$ is a non-negative even number, and it is zero otherwise. Here $b_{n,l,j}$ is a positive constant.

Summarizing,
we have
 $$\lambda_l(k) =  (-1)^l \sum_{j=0\atop j\, \mathrm{ even}}^k {{ (m-2)/2} \choose  j/2} \, j!\,  c_{n,l,j}\,  b_{n,l,j}. $$
Note that inside the summation all the terms are non-negative, and there is at least one term that is strictly positive, since there is  $j$ that satisfies $l-k\le j \le m-2$.
Thus, we have proved  that  $\lambda_l(k) \ne 0$ if and only if $ k-l$ is even and $m\ge k-l+2$.

Note that we extended $\lambda_l(q)$ to all non-negative values of $q$, but we did not do this for $I_q f $. For $I_q f$ this is also possible, but we would need $f$ to be infinitely smooth (which we want to avoid). Instead of extending $I_q f$, we will proceed as follows.

Let $H_l$ be a spherical harmonic of order $l$. First consider the case $-1<\Re q<0$. Then
\begin{align*} \int_{S^{n-1}} A_{K,m, \xi}^{(q)} (0) H_l(\xi) \, d\xi &= \frac{m}{(m-q)(n-1)} \int_{S^{n-1}} I_q  \left( \rho_K^{m-q} \right) (\xi) H_l(\xi) \, d\xi \\
&=\frac{m}{(m-q)(n-1)} \int_{S^{n-1}} \rho_K^{m-q}(\xi) I_q \left(  H_l \right) (\xi) \, d\xi\\
& = \frac{m}{(m-q)(n-1)}  \lambda_l(q) \int_{S^{n-1}}   \rho_K^{m-q}   (\xi) H_l(\xi) \, d\xi.
\end{align*}
The left and right-hand  sides of the equality are analytic functions of $q$ in the domain $\{q\in \mathbb C: \Re(q)>-1 \mbox{ and } q\ne m\}$.
Thus they are equal when $q$ is an integer  $k>\max\{m,N\}$. Since $A_{K,m, \xi}$ is a polynomial of $t$ for every $\xi$, we have    $$ A_{K,m, \xi}^{(k)} (0) = 0,$$ for all $\xi\in S^{n-1}$, and so
$$\lambda_l(k) \int_{S^{n-1}}   \rho_K^{m-k}   (\xi) H_l(\xi) \, d\xi = 0.$$

If $l\ge k-m+2$ and $k-l$ is even, then $\lambda_l(k)\ne 0$, and thus $\rho_K^{m-k}$ does not have spherical harmonics of order $l$, where $l=k$(mod 2), in its expansion.
Thus if $k$ is even, then $\rho_K^{m-k}+\rho_{-K}^{m-k}$ is a combination
of spherical harmonics of even degrees less than $k-m+2$. If $k$ is odd, then $\rho_K^{m-k}-\rho_{-K}^{m-k}$ is a combination
of spherical harmonics of odd degrees less than $k-m+2$.

Note that spherical harmonics are restrictions of homogeneous harmonic polynomials to the sphere.  Thus, as a function on $\mathbb R^n$, we
have, when $k$ is even,
$$ \|x\|_K^{k-m} + \| - x\|_K^{k-m} = \sum_{l=0\atop l\, \mathrm{ even} }^{k-m} H_l\left(\frac{x}{|x|}\right) |x|^{k-m}=\sum_{l=0}^{k-m} H_l\left(x\right) |x|^{k-m-l}.$$
Since $k-m-l$ is even, the right-hand side is a polynomial on $\mathbb R^n$.

Similarly, if $k$ is odd,
$$ \|x\|_K^{k-m} - \| - x\|_K^{k-m} = \sum_{l=0\atop l\, \mathrm{ odd} }^{k-m} H_l\left(\frac{x}{|x|}\right) |x|^{k-m}=\sum_{l=0}^{k-m} H_l\left(x\right) |x|^{k-m-l}$$
is also a polynomial on $\mathbb R^n$.

 Now we proceed  as in the proof of Theorem 3.7 in \cite{KMY}. It was shown there that if $K$ is a convex body such that $\|x\|_K^{l} + \| - x\|_K^{l}$ is a polynomial for even $l$ and $\|x\|_K^{l} - \| - x\|_K^{l}$ is a polynomial for odd $l$, then  $K$ is an ellipsoid.

\ep

\section{Final remarks}
One can also consider the  function $P_{K,m,\xi} (t)$ that gives the $m$th intrinsic volume of the section $K\cap \{\xi^\perp + t\xi\}$,
$$P_{K,m,\xi} (t) = V_m (K\cap \{\xi^\perp + t\xi\}).$$

If $m$ is even and $K$ is an ellipsoid, it is easy to see that $P_{K,m,\xi} (t)$ is a polynomial of $t$ of degree $m$ (on its support), for every $\xi$. The argument goes as in Lemma \ref{Lem}.
Using the Kubota integral recursion formula, we have
$$P_{K,m, \xi} (t) = \frac{\kappa_{n-1}}{\kappa_m}\int_{G( \xi^\perp, m)} V_m ((K| (H\vee\xi))\cap \{\xi^\perp + t\xi\})  \, dH.$$

We  now use the fact that $K| (H\vee\xi)$ is an ellipsoid, and therefore polynomially integrable in $H\vee\xi$ if $m$ is even (since $H\vee\xi$ is $(m+1)$-dimensional).

Question. Are ellipsoids the only convex bodies in $\mathbb R^n$ for which $P_{K,m, \xi} (t)$ is a polynomial of $t$ for every $\xi$?

\end{document}